\documentclass[12pt,a4paper]{article}
\usepackage{latexsym}
\newtheorem{lem}{Lemma}
\newtheorem{thm}{Theorem}
\newenvironment{pf}{\textbf{Proof\ }}{\hfill$\Box$\smallskip}
\newcommand{\sub}{\subseteq}
\newcommand{\cl}[1]{\left\lceil{#1}\right\rceil}
\newcommand{\fl}[1]{\left\lfloor{#1}\right\rfloor}
\newcommand{\la}{\lambda}
\newcommand{\si}{\sigma}
\newcommand{\cd}{\mathcal{D}}
\newcommand{\cp}{\mathcal{P}}
\newcommand{\N}{\mathbf{N}}
\title{Combinatorial proofs of $q$-series identities}
\author{Robin Chapman\\
Division of Information and Communication Sciences\\ Macquarie University\\
North Ryde, NSW 2109, Australia\\ \texttt{rjc@maths.ex.ac.uk}}
\date{27 August 2001}
\begin{document}
\maketitle
\footnotetext{Present address: School of Mathematical Sciences, University
of Exeter, Exeter, EX4 4QE, UK}

\section{Introduction}

We provide combinatorial proofs of six of the ten $q$-series identities listed
in \cite[Theorem 3]{AJO}. Andrews, Jim\'enez-Urroz and Ono prove these
identities using formal manipulation of identities arising in the theory of
basic hypergeometric series. Our proofs are purely combinatorial,
based on interpreting
both sides of the identities as generating functions for certain partitions.
One of these identities arose in the work of Zagier \cite{Z} and the author
has already given a combinatorial proof of the same~\cite{C}.

\section{Notation}

For convenience we summarize the notations we use for partitions and sets
of partitions.

The generic partition is denoted by $\la$. The number
partitioned by $\la$ is~$N_\la$. The parts of $\la$ are
$\la_1\ge\la_2\ge\cdots$.
In particular the largest part of $\la$ is $\la_1$. We let $n_\la$
denote the number of parts in $\la$. Also we let $n_\la(d)$ be the number of
parts
in $\la$ which equal $d$, so that $n_\la=\sum_d n_\la(d)$. We let $d_\la$
denote
the number of distinct parts in~$\la$.

The set of all partitions is denoted by~$\cp$. The set of all partitions
into distinct parts is denoted by~$\cd$.

\section{The identities}

Each of the identities in Theorem~3 of \cite{AJO} has the following structure:
\begin{equation}\label{main}
\sum_{N=0}^\infty\left[\prod_{j=1}^\infty b_j(q)-
\prod_{j=1}^N b_j(q)\right]=\prod_{j=1}^\infty b_j(q)\sum_{d=1}^\infty c_d(q)
+G(q).
\end{equation}
Here the $b_j(q)$ are power series in $q$ tending $q$-adically to 1 as $j\to\infty$,
the $c_d(q)$ are power series in $q$ tending $q$-adically to 0 as $d\to\infty$
and $G(q)$ is an explicitly given power series. Our strategy is to interpret
$$\sum_{N=0}^\infty\left[\prod_{j=1}^\infty b_j(q)-
\prod_{j=1}^N b_j(q)\right]\qquad\mathrm{and}\qquad
\prod_{j=1}^\infty b_j(q)\sum_{d=1}^\infty c_d(q)$$
as weighted generating functions for two classes of partitions, and to give a
weight preserving ``almost-bijection'' beween the two classes. There may
be exceptional partitions which are not paired under the almost-bijection;
these account for the series $G(q)$.

\begin{thm}\label{mainthm}
The identity (\ref{main}) is valid for the following values of
$b_j(q)$, $c_d(q)$ and $G(q)$:

\begin{center}
\begin{tabular}{|c||c|c|c|}
\hline
Case&$b_j(q)$&$c_d(q)$&$G(q)$\\
\hline(i)&$1/(1-q^j)$&$q^d/(1-q^d)$&$0$\\
(ii)&$(1+q^j)/(1-q^j)$&$2q^d/(1-q^{2d})$&$0$\\
(iii)&$(1-q^{2j-1})/(1-q^{2j})$&$(-1)^dq^d/(1-q^d)$&$0$\\
(iv)&$1-q^j$&$q^d/(1-q^d)$&$
\sum_{r=1}^\infty(-1)^r[(3r-1)q^{r(3r-1)}+3rq^{r(3r+1)}]$\\
(v)&$(1-q^j)/(1+q^j)$&$2q^d/(1-q^{2d})$&$4\sum_{r=1}^\infty(-1)^rrq^{r^2}$\\
(vi)&$(1-q^{2j})/(1-q^{2j+1})$&$(-1)^dq^d/(1-q^d)$&$(1-q)\sum_{r=1}^\infty rq^{r(r+1)/2}$\\
\hline
\end{tabular}.
\end{center}
\end{thm}
These correspond to parts 1--3 and 8--10 of \cite[Theorem 3]{AJO}. The remaining
parts 4--7 have a mock theta function as $G(q)$ and they appear to lie beyond the
methods of this paper.

To prove the theorem  in each case we start with a rearrangement of
the left side of (\ref{main}).

\begin{lem}\label{sum}
Let $a_n(q)$ $(n=1,2,\ldots)$
be power series in the indeterminate $q$ with $a_n\to 0$ in the $q$-adic
topology. Then
\begin{equation}\label{sumlemm}
\sum_{N=0}^\infty\left[\prod_{j=1}^\infty(1+a_j(q))-\prod_{j=1}^N(1+a_j(q))\right]
=\sum_{n=1}^\infty na_n(q)\prod_{j=1}^{n-1}(1+a_j(q)).
\end{equation}
\end{lem}
\begin{pf}
The product
$$\prod_{j=1}^\infty(1+a_j(q))$$
is the sum of all terms $a_S(q)=\prod_{j\in S}a_j(q)$ where $S$ runs through
the finite subsets of $\N=\{1,2,\ldots\}$.
Similarly
$$\prod_{j=1}^N(1+a_j(q))$$
is the sum of all terms $a_S(q)$ where $S$ runs through
the finite subsets of $\{1,2,\ldots,N\}$. Hence
the difference
$$\prod_{j=1}^\infty(1+a_j(q))-\prod_{j=1}^N(1+a_j(q))$$
is the sum of all $a_S(q)$ over all finite $S$ having an element strictly larger
than~$N$. Each $a_S(q)$ occurs in this difference for exactly $\max(S)$
distinct~$N$, where $\max(S)$ is the largest element of~$S$. Hence
$$\sum_{N=0}^\infty\left[\prod_{j=1}^\infty(1+a_j(q))-\prod_{j=1}^N(1+a_j(q))\right]
=\sum_{S\sub\N\atop 0<|S|<\infty}\max(S)a_S(q).$$
But
$$a_n(q)\prod_{j=1}^{n-1}(1+a_j(q))=\sum_{S\sub\N\atop\max(S)=n}a_S(q)$$
and so (\ref{sumlemm}) follows.
\end{pf}

\noindent\textbf{Proof of Theorem~\ref{mainthm}\ }

\noindent Case (i)
By Lemma~\ref{sum}, (\ref{main}) is equivalent to
$$\sum_{n=1}^\infty\frac{nq^n}{1-q^n}\prod_{j=1}^{n-1}\frac{1}{1-q^j}
=\prod_{j=1}^\infty\frac{1}{1-q^j}\sum_{d=1}^\infty\frac{q^d}{1-q^d}.\eqno(3)$$
The coefficient of $q^m$ in
$$\frac{q^n}{1-q^n}\prod_{j=1}^{n-1}\frac{1}{1-q^j}$$
is the number of partitions of $m$ with largest part~$n$. Thus
$$\sum_{n=1}^\infty\frac{nq^n}{1-q^n}\prod_{j=1}^{n-1}\frac{1}{1-q^j}
=\sum_{\la\in\cp}\la_1q^{N_\la}.$$
Now
$$\frac{1}{1-q^d}\frac{q^d}{1-q^d}=\sum_{m=1}^\infty m q^{dm}$$
and so
$$\frac{q^d}{1-q^d}\prod_{j=1}^\infty\frac{1}{1-q^j}
=\sum_{\la\in\cp}n_\la(d)q^{N_\la}.$$
Hence
$$\prod_{j=1}^\infty\frac{1}{1-q^j}\sum_{d=1}^\infty\frac{q^d}{1-q^d}
=\sum_{\la\in\cp}n_\la q^{N_\la}.$$
As the number of parts in a partition equals the largest part of its
conjugate, then
$$\sum_{\la\in\cp}n_\la q^{N_\la}
=\sum_{\la\in\cp}\la_1 q^{N_\la}$$
which proves~(3).

\medskip
\noindent Case (ii)
In this case, using Lemma~\ref{sum}, (\ref{main}) becomes
$$\sum_{n=1}^\infty\frac{2nq^n}{1-q^n}\prod_{j=1}^{n-1}\frac{1+q^j}{1-q^j}
=2\prod_{j=1}^\infty\frac{1+q^j}{1-q^j}\sum_{d=1}^\infty\frac{q^d}{1-q^{2d}}.
\eqno(4)$$
Note that
$$\frac{1+q^j}{1-q^j}=1+2\sum_{m=1}^\infty q^{mj}.$$
Then
$$\frac{2q^n}{1-q^n}\prod_{j=1}^{n-1}\frac{1+q^j}{1-q^j}
=\sum_{\la\in\cp\atop \la_1=n}2^{d_\la}q^{N_\la}.$$
Hence
$$\sum_{n=1}^\infty\frac{2nq^n}{1-q^n}\prod_{j=1}^{n-1}\frac{1+q^j}{1-q^j}
=\sum_{\la\in\cp}\la_12^{d_\la}q^{n_\la}.$$
On the other side note that
$$\frac{1+q^d}{1-q^d}\frac{q^d}{1-q^{2d}}=\frac{q^d}{(1-q^d)^2}
=\sum_{m=1}^\infty m q^{md}.$$
It follows that
$$\frac{2q^d}{1-q^{2d}}\prod_{j=1}^\infty\frac{1+q^j}{1-q^j}
=\sum_{\la\in\cp}n_\la(d)2^{d_\la}q^{N_\la}$$
and so
$$2\prod_{j=1}^\infty\frac{1+q^j}{1-q^j}\sum_{d=1}^\infty\frac{q^d}{1-q^{2d}}
=\sum_{\la\in\cp}n_\la2^{d_\la}q^{N_\la}.$$
As the number of distinct parts in a permutation is invariant under
conjugation, (4) follows. (Note that in the Ferrers diagram of a partition,
the number of distinct parts is the number of \emph{extreme} boxes---those
boxes which are both at the end of a row and the bottom of a column. Their
number is manifestly invariant under conjugation).

\medskip
\noindent Case (iii)
By Lemma~\ref{sum}, (\ref{main}) is equivalent in this case to
$$\sum_{n=1}^\infty n\frac{-q^{2n-1}+q^{2n}}{1-q^{2n}}
\prod_{j=1}^{n-1}\frac{1-q^{2j-1}}{1-q^{2j}}
=\prod_{j=1}^\infty\frac{1-q^{2j-1}}{1-q^{2j}}
\sum_{d=1}^\infty(-1)^d\frac{q^d}{1-q^d}.\eqno(5)$$
Let $\cp^o$ denote the set of partitions with no repeated odd parts, and
let $n^o_\la$ denote the number of odd parts of~$\la$. Then
$$\prod_{j=1}^\infty\frac{1-q^{2j-1}}{1-q^{2j}}
=\sum_{\la\in\cp^o}(-1)^{n^o_\la}q^{N_\la}$$
while
$$\prod_{j=1}^n\frac{1-q^{2j-1}}{1-q^{2j}}
=\sum_{\la\in\cp^o\atop \la_1\le 2n}(-1)^{n^o_\la}q^{N_\la}.$$
It follows that
$$\frac{-q^{2n-1}+q^{2n}}{1-q^{2n}}
\prod_{j=1}^{n-1}\frac{1-q^{2j-1}}{1-q^{2j}}
\sum_{\la\in\cp^o\atop \la_1\in\{2n-1,2n\}}(-1)^{n^o_\la}q^{N_\la}.$$
Hence
$$\sum_{n=1}^\infty n\frac{-q^{2n-1}+q^{2n}}{1-q^{2n}}
\prod_{j=1}^{n-1}\frac{1-q^{2j-1}}{1-q^{2j}}
=\sum_{\la\in\cp^o}\cl{\la_1/2}(-1)^{n^o_\la}q^{N_\la}.$$
Let us consider the other side. For odd~$d$,
$$\frac{(-1)^dq^d}{1-q^d}(1-q^d)=-q^d$$
and for even~$d$,
$$\frac{(-1)^dq^d}{1-q^d}\frac{1}{1-q^d}=\frac{q^d}{(1-q^d)^2}
=\sum_{m=1}^\infty mq^{md}.$$
In all cases then
$$\frac{(-1)^dq^d}{1-q^d}\prod_{j=1}^\infty\frac{1-q^{2j-1}}{1-q^{2j}}
=\sum_{\la\in\cp^o}n_\la(d)(-1)^{n^o_\la}q^{N_\la}$$
and so
$$\prod_{j=1}^\infty\frac{1-q^{2j-1}}{1-q^{2j}}
\sum_{d=1}^\infty\frac{(-1)^dq^d}{1-q^d}
=\sum_{\la\in\cp^o}n_\la(-1)^{n^o_\la}q^{N_\la}.$$
The desired identity will follow when we construct an involution
$\si$ on $\cp^o$ with the properties that $\si$ preserves
$N_\la$ while interchanging $n_\la$ and $\cl{\la_1/2}$.

To construct $\si$ we construct a \emph{diagram} for each
$\la\in\cp^o$. As an example let $\la=(8,7,5,4,4,3,2,2,2,1)$.
Each part of $\la$ will yield a row in the diagram. An even part
$2k$ will give a row of $k$ 2s, while an odd part $2k+1$ will give
a row of $k$ 2s followed by a 1. This particular $\la$ gives the
diagram
$$\begin{array}{cccc}
2&2&2&2\\
2&2&2&1\\
2&2&1&\\
2&2&&\\
2&2&&\\
2&&&\\
2&&&\\
2&&&\\
1&&&
\end{array}.$$
The diagram determines the partition via its row sums. It is a Ferrers
diagram with the boxes replaced by 2s and 1s. The 1s can only occur
at the end of rows, and for $\la\in\cp^o$ they can only occur at the
bottom of columns since no odd part is repeated. Any such diagram of
2s and 1s with the 1s occurring only in extreme boxes comes from a partition
in~$\cp^o_\la$. We identify elements of $\cp^o_\la$ with their diagrams.
Now we define $\si$ to be conjugation of diagrams. The number of rows
in the diagram of $\la$ is $n_\la$ while the number of columns is
$\cl{\la_1/2}$. Thus $\si$ has the required properties and (5)
is proved.

\medskip
\noindent Case (iv)
This proof appeared in \cite{C}. We outline it for completeness.
In this case, using Lemma~\ref{sum}, (\ref{main})
becomes
\begin{eqnarray*}
&&\sum_{n=1}^\infty-nq^n\prod_{j=1}^{n-1}(1-q^j)
-\prod_{j=1}^\infty(1-q^j)\sum_{d=1}^\infty\frac{q^d}{1-q^d}\\
&=&\sum_{r=1}^\infty[(3r-1)q^{r(3r-1)/2}+3rq^{r(3r+1)/2}].
\qquad\qquad(6)
\end{eqnarray*}
Then
$$-q^n\prod_{j=1}^{n-1}(1-q^j)=\sum_{\la\in\cd\atop\la_1=n}
(-1)^{n_\la}q^{N_\la}$$
and so
$$\sum_{n=1}^\infty-nq^n\prod_{j=1}^{n-1}(1-q^j)
=\sum_{\la\in\cd}\la_1(-1)^{n_\la}q^{N_\la}.$$
Also as
$$-\frac{q^d}{1-q^d}(1-q^d)=-q^d$$
then
$$\frac{-q^d}{1-q^d}\prod_{j=1}^\infty(1-q^j)
=\sum_{\la\in\cd}n_\la(d)(-1)^{n_\la}q^{N_\la}$$
and so
$$-\prod_{j=1}^\infty(1-q^j)\sum_{d=1}^\infty\frac{q^d}{1-q^d}
=\sum_{\la\in\cd}n_\la(-1)^{n_\la}q^{N_\la}.$$
Hence
$$\sum_{n=1}^\infty-nq^n\prod_{j=1}^{n-1}(1-q^j)
-\prod_{j=1}^\infty(1-q^j)\sum_{d=1}^\infty\frac{q^d}{1-q^d}
=\sum_{\la\in\cd}(\la_1+n_\la)(-1)^{n_\la}q^{N_\la}.$$

We now apply Franklin's involution (\cite{F}). For this purpose a partition
$\la\in\cd$ is \emph{exceptional} if either it is empty, or has the form
$(2r-1,2r-2,\ldots,r+1,r)$ or $(2r,2r-1,\ldots,r+2,r+1)$ for some
positive integer~$r$. Franklin's involution $\si$ is defined
on the nonexceptional partitions in~$\cd$. For $\la\in\cd$ let $s_\la$
be its smallest part, and $t_\la$ be the largest integer $t$ with
$\la_t=\la_1+1-t>0$. Define $\si(\la)$ for non-exceptional $\la\in\cd$
as follows: if $s_\la\le t_\la$ remove the smallest part from $\la$
and add 1 to each of its $s_\la$ largest parts, while if
$s_\la>t_\la$ subtract one from each of its
$t_\la$ largest parts and create a new part~$t_\la$. The map $\si$ is
an involution. It negates $(-1)^{n_\la}$ and it preserves $N_\la$
and $\la_1+n_\la$. The sum of $(-1)^{n_\la}q^{N_\la}$
(respectively $(\la_1+n_\la)(-1)^{n_\la}q^{N_\la}$) over
all $\la\in\cd$ thus equals the sum over the exceptional~$\la$.
We recover the pentagonal number identity
$$\prod_{j=1}^\infty(1-q^j)=\sum_{\la\in\cd}(-1)^{n_\la}q^{n_\la}
=1+\sum_{r=1}^\infty(-1)^r[q^{r(3r-1)/2}+q^{r(3r+1)/2}]$$
and the identity
$$\sum_{\la\in\cd}(\la_1+n_\la)(-1)^{n_\la}q^{n_\la}
=\sum_{r=1}^\infty(-1)^r[(3r-1)q^{r(3r-1)/2}+3rq^{r(3r+1)/2}]$$
which we have seen is equivalent to (6).

\medskip
\noindent Case (v)
By Lemma~\ref{sum} (\ref{main}) is equivalent to
$$\sum_{n=1}^\infty\frac{-2nq^n}{1+q^n}
\prod_{j=1}^{n-1}\frac{1-q^j}{1+q^j}
-2\prod_{j=1}^\infty\frac{1-q^j}{1+q^j}
\sum_{d=1}^\infty\frac{q^d}{1-q^{2d}}
=4\sum_{r=1}^\infty (-1)^rrq^{r^2}.\eqno(7)$$
Now
$$\frac{1-q^j}{1+q^j}=1+2\sum_{m=1}^\infty(-1)^mq^{mj}$$
and
$$\frac{-2q^n}{1+q^n}=2\sum_{m=1}^\infty(-1)^mq^{mn}.$$
Hence
$$\frac{-2q^n}{1+q^n}\prod_{j=1}^{n-1}\frac{1-q^j}{1+q^j}
=\sum_{\la\in\cp\atop\la_1=n}(-1)^{n_\la}2^{d_\la}q^{N_\la}$$
and so
$$\sum_{n=1}^\infty\frac{-2nq^n}{1+q^n}
\prod_{j=1}^{n-1}\frac{1-q^j}{1+q^j}
=\sum_{\la\in\cp}\la_1(-1)^{n_\la}2^{d_\la}q^{N_\la}.$$
On the other hand
$$\frac{-2q^d}{1-q^{2d}}\frac{1-q^d}{1+q^d}
=\frac{-2q^d}{(1+q^d)^2}=2\sum_{m=1}^\infty(-1)^mmq^{md}$$
so that
$$\frac{-2q^d}{1-q^{2d}}\prod_{j=1}^\infty\frac{1-q^j}{1+q^j}
=\sum_{\la\in\cp}n_\la(d)(-1)^{n_\la}2^{d_\la}q^{N_\la}$$
and then
$$-2\prod_{j=1}^\infty\frac{1-q^j}{1+q^j}\sum_{d=1}^\infty\frac{q^d}{1-q^{2d}}
=\sum_{\la\in\cp}n_\la(-1)^{n_\la}2^{d_\la}q^{N_\la}.$$
We need to prove that
$$\sum_{\la\in\cp}(\la_1+n_\la)(-1)^{n_\la}2^{d_\la}q^{N_\la}
=4\sum_{r=1}^\infty(-1)^rrq^{r^2}.$$
Declare a partition \emph{exceptional} if it is either empty or
is ``square'', consisting of exactly $r$ parts all equal to~$r$.
The sum of $(\la_1+n_\la)(-1)^{n_\la}2^{d_\la}q^{N_\la}$ over the
exceptional $\la$ is $4\sum_{r=1}^\infty(-1)^rrq^{r^2}$, so all
we need to show is that the sum over the nonexceptional $\la$
vanishes. This time there is no convenient involution, so we need
a more subtle approach.

Let $s_\la$ be the smallest part of the partition $\la$
and let $s'_\la$ be the smallest part of its conjugate. In other words
$s'_\la$ is the number of parts in $\la$ equal to~$\la_1$. Suppose
that $\la$ is not exceptional. If $s_\la\le s'_\la$ define the \emph{right
neighbour} of $\la$ as the partition formed by removing the smallest part
of $\la$ and adding 1 to its $s_\la$ largest parts. If $s_\la\ge s'_\la$
define the \emph{left neighbour} of $\la$ as the partition formed by
subtracting 1 from its $s'_\la$ largest parts and creating a new part
$s'_\la$. If $\la'$ is the left neighbour of $\la$ then $\la$ is the
right neighbour of $\la'$ and conversely. Each nonexceptional $\la$ has
one or two neighbours. If we define a graph $G$ with vertex set the
non-exceptional partitions and joining each partition to its neighbours
then each vertex in $G$ has degree 1 or 2. Also $G$ is acyclic since
moving to a partition's right neighbour decreases $n_\la$ by 1.
Thus the graph $G$ is a disjoint union of paths of length at least 2.
For example the following partitions form a path moving from left to right:
(4,4,4,3,2,2,2), (5,5,4,3,2,2), (6,6,4,3,2) and (7,7,4,3).
In each path $N_\la$ and $\la_1+n_\la$ are constant while the $(-1)^{n_\la}$
alternate in sign. We also see that if one endpoint has $d$ distinct parts,
then all non-endpoints have $d+1$ distinct parts and the other
endpoint has $d$ distinct parts. Thus the sums of $(-1)^{n_\la}2^{d_\la}$
and $(\la_1+n_\la)(-1)^{n_\la}2^{d_\la}$ over each path vanish. Hence
the sums of $(-1)^{n_\la}2^{d_\la}q^{N_\la}$
and $(\la_1+n_\la)(-1)^{n_\la}2^{d_\la}q^{N_\la}$ over the nonexceptional
partitions vanish. We thus recover the classical theta function
factorization
$$\prod_{j=1}^\infty\frac{1-q^j}{1+q^j}=
\sum_{\la\in\cp}(-1)^{n_\la}2^{d_\la}q^{N_\la}
=1+2\sum_{r=1}^\infty(-1)^rq^{r^2}.$$
the identity
$$\sum_{\la\in\cp}(\la_1+n_\la)(-1)^{n_\la}2^{d_\la}q^{N_\la}
=4\sum_{r=1}^\infty(-1)^rrq^{r^2}$$
which is equivalent to (7).

A different combinatorial proof of this theta function factorization
appears in~\cite{A2}.

\medskip
\noindent Case (vi)
In this last case, applying Lemma~\ref{sum}, (\ref{main}) becomes
\begin{eqnarray*}
\frac{1}{1-q}
\sum_{n=1}^{\infty}n\frac{-q^{2n}+q^{2n+1}}{1-q^{2n+1}}
\prod_{j=1}^{n-1}\frac{1-q^{2j}}{1-q^{2j+1}}
&-&\prod_{j=1}^\infty\frac{1-q^{2j}}{1-q^{2j-1}}
\sum_{d=1}^\infty(-1)^d\frac{q^d}{1-q^d}\\
&=&\sum_{r=1}^\infty rq^{r(r+1)/2}.\qquad\qquad(8)
\end{eqnarray*}
Let $\cp^e$ denote the set of partitions with no repeated even part
and $n^e_\la$ the number of even parts in the partition $\la$. Then
$$\frac{1}{1-q}\prod_{j=1}^\infty\frac{1-q^{2j}}{1-q^{2j+1}}
=\sum_{\la\in\cp^e}(-1)^{n^e_\la}q^{N_\la},$$
$$\frac{1}{1-q}\prod_{j=1}^n\frac{1-q^{2j}}{1-q^{2j+1}}
=\sum_{\la\in\cp^e\atop\la_1\le 2n+1}(-1)^{n^e_\la}q^{N_\la}$$
and
$$\frac{1}{1-q}\frac{-q^{2n}+q^{2n+1}}{1-q^{2n+1}}
\prod_{j=1}^{n-1}\frac{1-q^{2j}}{1-q^{2j+1}}
=\sum_{\la\in\cp^e\atop\la_1\in\{2n,2n+1\}}(-1)^{n^e_\la}q^{N_\la}.$$
Consequently
$$\frac{1}{1-q}\sum_{n=1}^\infty n\frac{-q^{2n}+q^{2n+1}}{1-q^{2n+1}}
\prod_{j=1}^{n-1}\frac{1-q^{2j}}{1-q^{2j+1}}
=\sum_{\la\in\cp^e}\fl{\la_1/2}(-1)^{n^e_\la}q^{N_\la}.$$
On the other hand for $d$ even
$$-\frac{(-1)^dq^d}{1-q^d}(1-q^d)=-q^d$$
while for $d$ odd
$$-\frac{(-1)^dq^d}{1-q^d}\frac{1}{1-q^d}=\frac{q^d}{(1-q^d)^2}
=\sum_{m=1}^\infty mq^{md}.$$
Hence
$$-(-1)^d\frac{q^d}{1-q^d}\prod_{j=1}^\infty\frac{1-q^{2j}}{1-q^{2j-1}}
=\sum_{\la\in\cp^e}n_\la(d)(-1)^{d^e_\la}q^{N_\la}$$
and so
$$-\prod_{j=1}^\infty\frac{1-q^{2j}}{1-q^{2j-1}}
\sum_{d=1}^\infty(-1)^d\frac{q^d}{1-q^d}
=\sum_{\la\in\cp^e}n_\la(-1)^{d^e_\la}q^{N_\la}.$$
We thus need to prove that
$$\sum_{\la\in\cp^e}(\fl{\la_1/2}+n_\la)(-1)^{d^e_\la}q^{N_\la}
=\sum_{r=1}^\infty rq^{r(r+1)/2}.$$
We define $\la\in\cp^e$ to be \emph{exceptional} if it is empty
or it has the form $((2s-1)^s)$ or $((2s+1)^s)$ for some positive
integer~$s$. The sum of
$(\fl{\la_1/2}+n_\la)(-1)^{d^e_\la}q^{N_\la}$ over the exceptional
$\la$ is $\sum_{r=1}^\infty rq^{r(r+1)/2}$. We must show that the
corresponding sum over the nonexceptional $\la$ vanishes, and for
this we need to define a suitable involution $\si$ on the set
of nonexceptional $\la\in\cp^e$.

We represent each $\la\in\cp^e$ by a \emph{diagram}. The parts of
$\la$ in descending order determine the rows of the diagram in
descending order. An odd part $2k+1$ becomes a row consisting
of 1 followed by $k$ 2s. An even part $2k$ becomes a row consisting
of 1 followed by $k-1$ 2s and ending with a 1. For instance
$(9,8,5,5,5,4,3,2,1,1)$ becomes
$$\begin{array}{ccccc}
1&2&2&2&2\\
1&2&2&2&1\\
1&2&2&&\\
1&2&2&&\\
1&2&2&&\\
1&2&1&&\\
1&2&&&\\
1&1&&&\\
1&&&&\\
1&&&&
\end{array}.$$
A diagram is a Ferrers diagram filled with 1s and 2s. The left column
is all 1s while the remainder of the diagram is filled with 2s apart
from the extreme boxes, which are filled with either 1s or 2s.
(We cannot have a 1 in any column but the first except at the bottom
as there are no repeated even parts). The Ferrers diagram has
$n_\la$ rows and $1+\fl{\la_1/2}$ columns. Let $s_\la$ be the smallest part
of the partition $\la$ and if $\la\in\cp^e$ let $s_\la''$ be the
sum of the final column in the associated diagram. We define $\si(\la)$
as follows: if $s_\la<s''_\la$ or if $s_\la=s''_\la$ is even
remove the last part of $\la$ and
add a new final column with sum $s_\la$ to the diagram of $\la$.
We then get the diagram of the partition $\si(\la)$. If
$s_\la>s''_\la$ or if $s_\la=s''_\la$ is odd
remove the last column of the diagram of $\la$
and create a new final row with sum~$s''_\la$. Again,
this is the diagram of the partition $\si(\la)$. We can check that $\si$
is an involution preserving $N_\la$ and negating $(-1)^{n_\la}$.
As $\si$ preserves the sum of the number of rows and the number of columns
in the diagram it also preserves $\fl{\la_1/2}+n_\la$. We thus recover the
classical theta function factorization
$$\prod_{j=1}^\infty\frac{1-q^{2j}}{1-q^{2j-1}}
=\sum_{\la\in\cp^e}(-1)^{n^e_\la}q^{N_\la}
=\sum_{r=0}^\infty q^{r(r+1)/2}$$
and our desired identity
$$\sum_{\la\in\cp^e}(\fl{\la_1/2}+n_\la)(-1)^{n^e_\la}q^{N_\la}
=\sum_{r=1}^\infty rq^{r(r+1)/2}$$
which is equivalent to (8).

This proof of the theta function factorization is essentially equivalent to
that in~\cite{A2}.

\hfill$\Box$\smallskip

\section{Mock theta identities}

Four of the identities in \cite[Theorem 3]{AJO} involve mock theta
functions. The proofs of these identities appear to be beyond the scope
of the present methods, but these methods do provide combinatorial
interpretations of the identities.
We consider one case in detail: the fourth identity.
This identity is equivalent to the following identity from Ramanujan's
``Lost'' Notebook which was proved by Andrews \cite{A1}:
\begin{eqnarray*}
\sum_{N=0}^\infty\left[\prod_{j=1}^\infty(1+q^j)-\prod_{j=1}^N(1+q^j)\right]
&=&\prod_{j=1}^\infty(1+q^j)\left[-\frac12+\sum_{d=1}^\infty\frac{q^d}{1-q^d}
\right]\\
&&+\frac12\sum_{n=0}^\infty\prod_{j=1}^n\frac{q^j}{(1+q^j)}.\qquad\qquad(9)
\end{eqnarray*}
By Lemma~1
$$\sum_{N=0}^\infty\left[\prod_{j=1}^\infty(1+q^j)-\prod_{j=1}^N(1+q^j)\right]
=\sum_{n=1}^\infty nq^n\prod_{j=1}^{n-1}(1+q^j).$$
But
$$\sum_{n=1}^\infty q^n\prod_{j=1}^{n-1}(1+q^j)=
\sum_{\la\in\cd\atop \la_1=n}q^{N_\la}$$
and so
$$\sum_{n=1}^\infty nq^n\prod_{j=1}^{n-1}(1+q^j)
=\sum_{\la\in\cd}\la_1q^{N_\la}.$$
Next
$$\frac{q^d}{1-q^d}(1+q^d)=q^d+2\sum_{m=2}^\infty q^{md}.$$
It follows that
$$\frac{q^d}{1-q^d}\prod_{j=1}^\infty(1+q^j)
=\sum_{\la\in\cd_{d,1}}q^{N_\la}
+\sum_{\la\in\cd_{d,2+}}q^{N_\la}$$
where $\cd_{d,1}$ is the set of $\la\in\cd$ having $d$ as a part
while $\cd_{d,2+}$ is the set of partitions $\la$ having at least two
occurences of the part $d$ but with no other repeated parts. It follows that
$$\prod_{j=1}^\infty(1+q^j)\frac{q^d}{1-q^d}
=\sum_{\la\in\cd}n_\la q^{N_\la}+2\sum_{\la\in\cd_{(1)}}q^{N_\la}$$
where $\cd_{(1)}$ denotes the set of partitions having exactly one
part with multiplicity at least~2. Of course
$$\prod_{j=1}^\infty(1+q^j)=\sum_{\la\in\cd}q^{N_\la}.$$

For each $n$
$$\prod_{j=1}^n\frac{q^j}{(1+q^j)}
=\sum_{\la\in\cp_{(n)}}(-1)^{n_\la-n}q^{N_\la}$$
where $\cp_{(n)}$ is the set of partitions having no parts greater
than $n$ and having all numbers between 1 and $n$ as parts. These
partitions have a complicated description, but their conjugates
have a nicer one: $\la\in\cp_{(n)}$ if and only if its conjugate
lies in $\cd$ and has exactly $n$ parts. Hence
$$\sum_{n=0}^\infty\prod_{j=1}^n\frac{q^j}{(1+q^j)}
=\sum_{\la\in\cd}(-1)^{\la_1-n_\la}q^{N_\la}.$$
Putting all this together, the identity under consideration is equivalent
to
$$\sum_{\la\in\cd}
\left[\la_1-n_\la+\frac{1-(-1)^{\la_1-n_1}}{2}\right]q^{N_\la}
=2\sum_{\la\in\cd_{(1)}}q^{N_\la}$$
or equivalently
$$\sum_{\la\in\cd}\cl{(\la_1-n_\la)/2}q^{N_\la}
=\sum_{\la\in\cd_{(1)}}q^{N_\la}.$$
For a partition $\la$ the quantity $r_\la=\la_1-n_\la$ is often
called its \emph{rank}. The combinatorial interpretation of this identity
is that the sum of $\cl{r_\la/2}$ over all partitions $\la$ of $n$
into distinct parts equals the number of partitions of $n$ having exactly
one repeated part.
For instance here are the partitions $\la$ of 8 into distinct parts.
\begin{center}
\begin{tabular}{|c|c|c|}
\hline
$\la$&$r_\la$&$\cl{r_\la/2}$\\
\hline
8&7&4\\
71&5&3\\
62&4&2\\
53&3&2\\
521&2&1\\
431&1&1\\
\hline
\end{tabular}
\end{center}
The sum of the $\cl{r_\la/2}$ is 13, and there are 13 partitions of
2 in which exactly one part is repeated namely, $61^2$, $51^3$, $4^2$, $42^2$,
$421^2$, $41^4$, $3^22$, $32^21$, $321^3$, $31^5$, $2^4$, $21^6$ and $1^8$.
It would be interesting to have a bijective proof of this fact.

\section{Acknowledgment}

This work was completed when the author was a Visiting Fellow
in the Division of Information and Communication Sciences at Macquarie
University. He wishes to thank Alf van der Poorten for his invitation to
Macquarie
and also George Andrews and Don Zagier for providing copies of their papers.

\end{document}